\documentclass{amsproc}

\usepackage[all]{xy}
\usepackage{slashed}
\usepackage{hyperref}

\copyrightinfo{2023}{American Mathematical Society}

\newtheorem{theorem}{Theorem}[section]

\theoremstyle{definition}
\newtheorem{definition}[theorem]{Definition}
\newtheorem{example}[theorem]{Example}

\theoremstyle{remark}
\newtheorem{remark}[theorem]{Remark}

\numberwithin{equation}{section}

\begin{document}

\title{A stroll in equivariant $K$-theory}


%    author one information
%\author{}
%\address{}
%\curraddr{}
%\email{}
%\thanks{}

%    author two information
\author{Chi-Kwong Fok}
\address{NYU-ECNU Institute of Mathematical Sciences, New York University Shanghai, 
3663 Zhongshan Road North,
Shanghai 200062, China}
\curraddr{Xi'an Jiaotong-Liverpool University, 111 Ren'ai Road, Suzhou Industrial Park, Suzhou, Jiangsu Province 215123, China}
\email{ChiKwong.Fok@xjtlu.edu.cn}
\thanks{We would like to thank Loring W. Tu for organizing the Special Session on Equivariant Cohomology, from which we learnt a lot about the current developments of the subject, and which led to the writing of the present article. We are also grateful to the referee for their meticulous comments and suggestions which greatly improve the exposition of the article.}

%\subjclass[2020]{Primary: 19L47}
%    The 2010 edition of the Mathematics Subject Classification is
%    now available.  If you are citing a classification from the
%    new scheme, use the following input coding instead.
\subjclass[2020]{Primary 19L47}
\keywords{Equivariant $K$-theory, Segal localization theorem, Atiyah-Segal localization formula, equivariant formality, twisted $K$-theory}

\date{}

\begin{abstract}
Equivariant $K$-theory is a generalized equivariant cohomology theory which is a hybrid of the $K$-theory of a topological
space and the representation theory of the group acting on it. In this article, we review the basics of equivariant $K$-theory and focus on the localization theorem and formula due to Atiyah and Segal, which have become important tools in equivariant topology nowadays. We then discuss the application of equivariant $K$-theory to equivariant formality, and briefly mention some recent developments. 
\end{abstract}

\maketitle
\tableofcontents
\section{Introduction}
Equivariant $K$-theory, introduced by Segal in his seminal PhD thesis \cite{Se2}, occupies a prominent position among the machinery in the study of equivariant topology. A natural extension of $K$-theory, equivariant $K$-theory is defined using equivariant vector bundles and thus amalgamates ordinary $K$-theory with representation theory which keeps track of group actions. This article, which is an elaboration of the expository talk presented by the author in the Special Session on Equivariant Cohomology in the AMS Spring Eastern Virtual Sectional Meeting 2022, aims to journey through some important results and ideas of equivariant $K$-theory which have become classical nowadays, and a brief account of some recent relevant work. 

After starting from the ground up with preliminaries and examples of (equivariant) $K$-theory, we focus on the Segal localization theorem and the Atiyah-Segal localization formula (\cite{Se2, AtSe}), which extract global information (i.e. equivariant $K$-theory and indices of the whole space) through fixed point data (topological localization) and algebraic localization of the representation ring. We will see that the representation ring as the coefficient ring of equivariant $K$-theory has a richer algebraic structure than that of the ordinary $K$-theory which makes localization work. Equally crucial is the invertibility of the Euler class of the normal bundle of the fixed point set in the localized $K$-theory ring, which is peculiar in the equivariant setting and key to the localization formula. These salient ideas of localization in fact have permeated the study in equivariant topology and often recur in many later influential works. These include the localization theorem for equivariant complex bordism \cite[Theorem 3.1]{tD}), the localization formula for equivariant (de Rham) cohomology (\cite{AB2}), the topology of moduli spaces of Riemann surfaces (\cite{AB}), and the localization of the infinite dimensional path integrals of quantum field theory (\cite{W}), to name a few. To illustrate the power of localization, we will, among other things, reproduce the delightfully elegant, index theoretic proof of the Weyl character formula.

We also discuss how equivariant $K$-theory enters in characterizing equivariant formality, a desirable property of group actions which is originally defined in terms of equivariant cohomology and commonly occurring in equivariant topology (\cite{GKM}). We introduce the notion of rational $K$-theoretic equivariant formality (RKEF), which roughly speaking amounts to the existence of stable equivariant structures of any vector bundle and is arguably more natural than the cohomological definition of equivariant formality (\cite{F}). We apply RKEF to establish new proofs of equivariant formality of some well-known examples, e.g. actions with connected maximal rank isotropy subgroups. We also discuss equivariant formality of compact homogeneous spaces $G/K$ with the left translation action by $K$ using RKEF, and provide a representation theoretic characterization of it (\cite{CF}). The latter has an algebro-geometric interpretation which has its root in the idea of `affinizing' Lie groups in \cite{Se}, which in turn plays a part in the localization theorem as well as the `delocalized' approach to constructing equivariant $K$-theory (\cite{RK}). At the end of this article we briefly mention some recent developments in equivariant twisted $K$-theory, which also bear the mark of the ideas of localization and affinization.
\section{Preliminaries}

\subsection{$K$-theory}
As the predecessor of equivariant $K$-theory, (topological) $K$-theory was introduced by Atiyah-Hirzebruch (\cite{AH1, AH2}) in order to formulate a generalization of the Grothendieck-Riemann-Roch Theorem in the differentiable setting, and was used as the vehicle for topological indices defined by Atiyah-Singer for their renowned index theorem (\cite{AS}). A powerful machinery in algebraic topology, $K$-theory probes into the topology of a space by classifying the vector bundles it can support, which gives an idea of how much one can continuously deform the `linear algebra' over the space, and the natural algebraic structures arising from such fiberwise operations on vector bundles as direct sum and tensor product.
\begin{definition}\label{defk}
	\begin{enumerate}
		\item Let $X$ be a compact Hausdorff space and $\text{Vect}_\mathbb{C}(X)$ the set of isomorphism classes of finite rank complex vector bundles over $X$. The $K$-theory ring $K^0(X)$ is defined to be the Grothendieck group of $\text{Vect}_\mathbb{C}(X)$, i.e. the stable isomorphism classes of formal differences of elements in $\text{Vect}_\mathbb{C}(X)$, i.e.
		\[\{[E]-[F]|[E], [F]\in\text{Vect}_\mathbb{C}(X)\}\left/\begin{pmatrix}[E_1]-[F_1]\sim [E_2]-[F_2]\text{ if there exists}\\ W\in\text{Vect}_\mathbb{C}(X)\text{ such that}\\ E_1\oplus F_2\oplus W\cong E_2\oplus F_1\oplus W\end{pmatrix}\right. .\]
		Addition and multiplication of $K^0(X)$ are given by the direct sum and tensor product of vector bundles respectively.
		\item\label{compsupp} For $X$ with a base point, the reduced $K$-theory $\widetilde{K}^0(X)$ is defined to be the kernel of the map $K^0(X)\to K^0(\text{pt})$ given by restricting vector bundles over $X$ to the base point.\footnote{The isomorphism class of the group $\widetilde{K}^0(X)$ is independent of the choice of the base point. The same is true of the ring structure of $\widetilde{K}^0(X)$ (to be defined after Definition \ref{defk}) for $X$ connected. However, if $X$ is disconnected, the ring structure of $\widetilde{K}^0(X)$ does depend on the connected component which contains the base point.} For a locally compact Hausdorff space $X$, define $X^+$ to be the one-point compactification of $X$ if it is non-compact, and the disjoint union of itself and a point if it is compact. The base point of $X^+$ is defined to be the extra point added to $X$. We then define $K^0(X)$ to be $\widetilde{K}^0(X^+)$. Also define the negative degree piece $K^{-n}(X)$ to be $K^0(\mathbb{R}^n\times X)$. If $Y$ is a closed subspace of $X$, then define $K^{-n}(X, Y):=K^{-n}(X\setminus Y)$.
	\end{enumerate}
\end{definition}
One can further define the external product 
\[K^{-n}(X)\otimes K^{-m}(Y)\to K^{-n-m}(X\times Y)\] 
by tensoring pullback bundles induced by the projections of $X\times Y$ onto the two factors. Letting $X=Y$ and composing the external product with the restriction to the diagonal yields the product structure of the $\mathbb{Z}$-graded ring $\displaystyle\bigoplus_{n=0}^\infty K^{-n}(X)$. 
\begin{example}\label{kex}
	\begin{enumerate}
		\item\label{kpt} Let $X$ be a point. Then any complex vector bundle over $X$ is a complex vector space, whose isomorphism class is determined by its dimension. So $\text{Vect}_\mathbb{C}(\text{pt})\cong\mathbb{N}$ and its Grothendieck group is isomorphic to $\mathbb{Z}$. Taking into account the tensor product of vector spaces, we get $K^0(\text{pt})$ is isomorphic to the ring of integers. On the other hand, $K^{-1}(\text{pt})$, which by definition is $\widetilde{K}^0(S^1)$, is 0 as any complex vector bundle on $S^1$ is trivial by the clutching construction of vector bundles on spheres in general, to be explained below. 
		\item In general any vector bundle over a sphere $S^n$ can be constructed by piecing the trivial bundles over the two hemispheres along the equator by way of a clutching function which maps from the equator $S^{n-1}$ to the structure group $\text{U}(m)$, whose homotopy class determines the isomorphism class of the vector bundle (\cite[Proposition 1.11]{Ha}). By a series of homotopies reducing the clutching function to a simpler form, it can be shown that any vector bundle over $S^2$, which can be thought of as $\mathbb{CP}^1$, is the direct sum of tensor powers of the hyperplane line bundle $H$, and that $K^0(S^2)$ is isomorphic to the ring $\mathbb{Z}[H]/((H-1)^2)$. The reduced $K$-theory $\widetilde{K}^0(S^2)$ then is isomorphic to the subring $\mathbb{Z}\cdot(1-H)$, and $K^{-2}(\text{pt}):=\widetilde{K}^0(S^2)$ is isomorphic to $\mathbb{Z}$ with generator corresponding to $1-H$. The negative of this generator, which corresponds to $1-H^{-1}$ (as $1-H+1-H^{-1}=-H^{-1}(H-1)^2=0$), is called the Bott class $\beta$ (\cite{At2}).\footnote{Out of an index theoretic consideration, the generator $1-H^{-1}$ instead of the other generator $1-H$ is christened the Bott class (see Example \ref{topindH}).}
\end{enumerate}
\end{example}
The reduced $K$-theory $\widetilde{K}^0(S^n)$ turns out to depend on the parity of $n$ only: it is $\mathbb{Z}$ for $n$ even and 0 for $n$ odd. This is the renowned Bott periodicity, first formulated equivalently as the intriguing periodicity phenomenon of the homotopy groups of the infinite unitary group and its classifying space (cf. \cite{Bo}). It also takes the following more general form. 
\begin{theorem}[Bott periodicity]There is an isomorphism $K^{-n}(X)\to K^{-n-2}(X)$ through the external product with $\beta$.
\end{theorem}
See \cite{At} or \cite[\S 2.1]{Ha} for the aforementioned proof which analyses the clutching function. Bott periodicity enables us to define $K^n(X)$ for any integer $n$ by identifying $K$-theory groups of degrees with the same parity. Now it makes sense to recast the $K$-theory ring as follows.
\begin{definition}
	The $K$-theory ring $K^*(X)$ is defined as the $\mathbb{Z}_2$-graded ring $K^0(X)\oplus K^{1}(X)$.
\end{definition}
$K$-theory resembles cohomology in many ways in that both are contravariant functors from the category of spaces to that of graded commutative rings, and $K$-theory satisfies most of the Eilenberg-Steenrod axioms for cohomology theory. For example, $K$-theory is homotopy invariant, because two homotopic maps pull the same vector bundle back to isomorphic vector bundles.\footnote{To be more precise, $K$-theory is invariant under \emph{proper} homotopy, as it is a compactly supported cohomology theory due to Definition \ref{defk} (\ref{compsupp}). For instance, $K^0(\mathbb{R})$ is not isomorphic to $K^0(\text{pt})$, though $\mathbb{R}$ is contractible. That is because $\mathbb{R}$ cannot be contracted to a point by any proper homotopy. } The only axiom $K$-theory fails to satisfy is the dimension axiom, as $K^2(\text{pt})$ is not 0 while the degree is greater than the dimension of a point. This makes $K$-theory the first example of generalized cohomology theories. 
\subsection{Equivariant $K$-theory}
Just like cohomology has Borel equivariant cohomology as a generalization for the category of spaces with group actions, $K$-theory also admits a natural equivariant counterpart. Introduced by Segal in his seminal PhD thesis \cite{Se2}, equivariant $K$-theory is defined using equivariant vector bundles, where the group action projects onto that on the base space and takes a fiber to another linearly. Below we list more precise definitions about equivariant $K$-theory, which basically is Definition \ref{defk} with every instance of vector bundles replaced by equivariant vector bundles. Though a straightforward extension from ordinary $K$-theory as it may seem, the use of equivariant vector bundles in equivariant $K$-theory allows us, as we will see, to capture with tools from representation theory the underlying equivariant topology, which may not be gleaned just from ordinary $K$-theory. 
\begin{definition}
Let $G$ be a compact Lie group which acts on a locally compact Hausdorff space $X$.\footnote{Compactness of the Lie group is essential for the reduced equivariant $K$-theory defined more generally as the stable equivalence classes of $\text{Vect}^G_\mathbb{C}(X)$ (see the next footnote) to be an abelian group. By compactness, one can invoke a variant of Peter-Weyl theorem to show that for any equivariant vector bundle $E$, there exists a complex $G$-representation $V$ and another equivariant vector bundle $E'$ such that $E\oplus E'\cong \textbf{V}$ (cf. \cite[Prop. 2.4]{Se2}). Then $[E']$ serves as the inverse of $[E]$ in the reduced equivariant $K$-theory group thus defined.}	
\begin{enumerate}
			\item If $V$ is a complex representation of $G$, then we use the boldfaced letter $\textbf{V}$ to denote the $G$-vector bundle $X\times V$ where $G$ acts diagonally.
			\item Let $\text{Vect}^G_\mathbb{C}(X)$ be the set of isomorphism classes of finite rank complex $G$-vector bundles over $X$. If $X$ is compact, the equivariant $K$-theory ring $K_G^0(X)$ is defined to be the Grothendieck group of $\text{Vect}^G_{\mathbb{C}}(X)$, with the abelian group structure given by $\oplus$ and the product structure given by $\otimes$. If further $X$ has a $G$-fixed point, then define the reduced equivariant $K$-theory $\widetilde{K}_G^0(X)$ to be the kernel of the map $K_G^0(X)\to K_G^0(\text{pt})$ given by restricting equivariant vector bundles over $X$ to the $G$-fixed point.\footnote{For $X$ not necessarily with a $G$-fixed point, $\widetilde{K}^0_G(X)$ is defined to be the stable equivalence classes of $\text{Vect}^G_\mathbb{C}(X)$, where $[E_1]$ is said to be stably equivalent to $[E_2]$ if there exist complex $G$-representations $V$ and $W$ such that $E_1\oplus\textbf{V}\cong E_2\oplus \textbf{W}$. When $X$ does have a $G$-fixed point, this definition is equivalent to the above definition via the map $[E]\mapsto[E]-[\textbf{V}]$, where $V$ is the restriction of $E$ at the $G$-fixed point.}
			\item If $X$ is non-compact, define $K_G^0(X)$ to be $\widetilde{K}_G^0(X^+)$ where $G$ fixes the point of compactification. Also define the negative degree piece $K_G^{-n}(X)$ to be ${K}_G^0(\mathbb{R}^n\times X)$, where $G$ acts on $\mathbb{R}^n$ trivially. If $Y$ is a $G$-invariant closed subspace of $X$, then define $K^{-n}_G(X, Y)$ to be $K^{-n}_G(X\setminus Y)$.
	\end{enumerate}
\end{definition}
As before one can define the external and internal product in equivariant $K$-theory of various degrees. There is also an equivariant version of Bott periodicity $K_G^{-n}(X)\cong K_G^{-n-2}(X)$ (\cite[Proposition 3.5]{Se2}) through the external product with the equivariant Bott class $\beta_G\in K_G^{-2}(\text{pt})$ (see Example \ref{ekex} (\ref{ekcoef}) below). So again this allows us to define $K_G^*(X)$ as the $\mathbb{Z}_2$-graded ring $K_G^0(X)\oplus K_G^1(X)$. Equivariant $K$-theory is an example of generalized equivariant cohomology theories (see for example \cite[Prop. 2.3, 2.6]{Se2} for equivariant (proper) homotopy invariance and the existence of the long exact sequence for a pair of spaces).
\begin{example}\label{ekex}
	\begin{enumerate}
		\item\label{ekcoef} Any equivariant complex $G$-vector bundle over a point is nothing but a complex representation of $G$. Thus $K^0_G(\text{pt})$ is the complex representation ring $R(G)$. More generally, $K^{-n}_G(\text{pt})$ is by definition $\widetilde{K}^{0}_G(S^n)$, where $G$ acts trivially. By the next example, $\widetilde{K}_G^0(S^n)\cong R(G)\otimes\widetilde{K}^0(S^n)$, which is $R(G)$ for $n$ even and 0 for $n$ odd. The equivariant Bott class $\beta_G\in K_G^{-2}(\text{pt})$ corresponds to $1\otimes(1-H^{-1})\in R(G)\otimes \widetilde{K}^0(S^2)\cong \widetilde{K}_G^0(S^2)$. 
		\item\label{ekhomog} If $H$ is a Lie subgroup of $G$, then $K_G^*(G/H)\cong K_H^*(\text{pt})\cong R(H)$. This is because any complex $G$-equivariant vector bundle over $G/H$ is determined by the restricted $H$-action on the fiber over the identity coset $e_GH$.
		\item If $\varphi\colon H\to G$ is a Lie group homomorphism and $X$ is a $G$-space, then it is also an $H$-space by defining $h\cdot x:=\varphi(h)\cdot x$. This induces a map $\text{Vect}_\mathbb{C}^G(X)\to \text{Vect}_\mathbb{C}^H(X)$ and hence a homomorphism $K_G^*(X)\to K_H^*(X)$. In particular, if $H$ is the trivial group, the map $f_G\colon K_G^*(X)\to K^*(X)$ is the forgetful map forgetting the equivariant structure of vector bundles. If further $X$ is a point, then the forgetful map becomes the augmentation map $R(G)\to \mathbb{Z}$, which sends a representation to its dimension. The kernel of the augmentation map is called the augmentation ideal and denoted by~$I(G)$.
		\item\label{isotypical} If $G$ acts on $X$ trivially, then $K_G^*(X)$ is isomorphic to $R(G)\otimes K^*(X)$ through the fiberwise isotypical decomposition, which sends an equivariant vector bundle $V$ to $\displaystyle \bigoplus_{W\in\text{Irr}(G)}W\otimes\text{Hom}_G(\textbf{W}, V)$. 
		\item On the other extreme, if $G$ acts freely on $X$, then $K_G^*(X)$ is isomorphic to $K^*(X/G)$ because $\text{Vect}_\mathbb{C}^G(X)$ and $\text{Vect}_\mathbb{C}(X/G)$ are in bijective correspondence through taking the quotient of an equivariant vector bundle by the free $G$-action, which admits as the inverse the pullback map induced by the point-orbit projection. This is consistent with our expectation from the motivation behind the Borel construction of equivariant cohomology that equivariant cohomology of a space with a free group action should be the ordinary cohomology of the quotient.
		\item The equivariant $K$-theory ring admits a natural $R(G)$-algebra structure through the map
		\begin{align*}
			R(G)&\to K_G^*(X)\\
			V&\mapsto \textbf{V}.
		\end{align*}
	\end{enumerate}
\end{example}
%The incorporation of the equivariant structure into $K$-theory makes the coefficient ring of equivariant $K$-theory, which is a representation ring as we see in Example \ref{ekex} (\ref{ekcoef}), richer in algebraic structure than that of ordinary $K$-theory. This has the added and crucial benefit of performing localization in the enlarged coefficient ring, which is the key idea behind the localization theorem and formula we will focus in latter sections. The notion of localization has permeated the study in equivariant topology, and is a primary source of inspiration for various similar results in other (generalized) equivariant cohomology theories (e.g. the localization formula for equivariant (de Rham) cohomology \cite{AB} and the localization theorem for equivariant complex bordism \cite[Theorem 3.1]{tD}).

\section{Localization theorem for equivariant $K$-theory}
In this section we shall introduce the first important result in equivariant $K$-theory, namely the Segal localization theorem (\cite[Proposition 4.1]{Se2}). It basically says that on localizing the coefficient ring, the equivariant $K$-theory restricts isomorphically to that of the fixed point set. Thus localization allows the fixed point set, which in many cases is easier to handle than the whole space, to capture important information about the equivariant $K$-theory of the whole space. Let us first state the abelian version of the localization theorem where the Lie group is a compact torus. 
\begin{theorem}[Segal localization theorem, abelian version]\label{segallocab}
	Let $T$ be a compact torus acting on a locally compact Hausdorff space $X$. Then the restriction map
	\[i^*\colon K_T^*(X)\to K_T^*(X^T)\cong R(T)\otimes K^*(X^T)\]
	becomes an isomorphism after localizing the coefficient ring $R(T)$, which is a ring of Laurent polynomials, at the zero ideal.
\end{theorem}
An immediate consequence is that the rank of $K_T^*(X)$ as a $R(T)$-module is the rank of $K^*(X^T)$ as an abelian group. For example, consider $S^2$ with $S^1$ acting by rotation. Then the fixed point set consists of two points which are the points of intersection of the sphere with the axis of rotation. It follows that the rank of $K_{S^1}^*(S^2)$ is the rank of $\displaystyle K^*(\text{pt}\cup\text{pt})\cong\mathbb{Z}\oplus\mathbb{Z}$, which is 2. We will verify this in the next section by computing $K_{S^1}^*(S^2)$ explicitly.

 Note that $X^T$ is a closed $T$-invariant subspace. From the long exact sequence for the pair $(X, X^T)$
\[\cdots\longrightarrow K_T^n(X\setminus X^T)\longrightarrow K_T^n(X)\stackrel{i^*}{\longrightarrow}K_T^n(X^T)\longrightarrow\cdots\]
and exactness of taking localization, we see that to show $i^*$ is an isomorphism after localization, it suffices to prove $K_T^*(X\setminus X^T)_{(0)}=0$, which can be reduced to proving that $K_T^*(Y)_{(0)}=0$ for any compact $T$-subspace $Y$ of $X\setminus X^T$.\footnote{Equivariant $K$-theory is continuous with respect to any directed system of relatively compact open $G$-subspaces (\cite[Corollary 2.12]{Se2}), so it suffices to show that $K_T^*(U)_{(0)}=0$ for $U$ a relatively compact open $T$-subspace. Using the long exact sequence associated to the pair $(\overline{U}, \overline{U}\setminus U)$ of compact $T$-subspaces enables one to obtain the claimed reduction.} Now we shall follow the proof strategy which consists of applying the Mayer-Vietoris sequence to a cover of $Y$ by slices of orbits, whose equivariant $K$-theory is then shown to have null contribution after localization. In more details, we have that, for each $T$-orbit $O_i$ in $Y$, there exists a closed neighborhood $S_i$ which equivariantly retracts onto $O_i$ by the slice theorem. This induces a map $K_T^*(O_i)\to K_T^*(S_i)$, making $K_T^*(S_i)$ a module over $K_T^*(O_i)$. By assumption the stabilizer subgroup of $O_i$ is a proper subgroup $T_i$ of $T$. Thus we have
\begin{align*}
	K_T^*(O_i)&\cong K_T^*(T/T_i)\\
			&\cong R(T_i)\ (\text{by Example \ref{ekex} (\ref{ekhomog})})\\
			&\cong\mathbb{Z}[\Lambda^*_{T_i}]\text{, where }\Lambda^*_{T_i}\text{is the character group of }T_i.
\end{align*}
Here comes the key observation: $R(T_i)_{(0)}$ as a $R(T)_{(0)}$-module is 0 because there exists a nonzero character $\alpha\in \Lambda^*_T$ lying in the kernel of the restriction $\Lambda^*_T\to\Lambda^*_{T_i}$, and it annihilates $R(T_i)$. In this way $K_T^*(O_i)_{(0)}$ and hence $K_T^*(S_i)_{(0)}$ are zero as $R(T)_{(0)}$-modules. By compactness choose finitely many $S_j$'s which cover $Y$. Localizing the Mayer-Vietoris sequence with respect to this finite cover then yields $K_T^*(Y)_{(0)}=0$ as desired.

The Segal localization theorem in its full generality concerns the action by a general compact Lie group $G$ and localization with respect to an arbitrary prime ideal $\mathfrak{p}$ of $R(G)$. Like the abelian version, the general one is also an assertion of an isomorphism upon restriction of equivariant $K$-theory and localization at $\mathfrak{p}$. To get an idea of what orbits we should restrict to, we shall, taking cue from the previous paragraph on the proof of the abelian version, answer the question: if $H$ is a subgroup of $G$ (supposedly the stabilizer of a certain point), what is a necessary and sufficient condition for $R(H)_\mathfrak{p}$ to be a nonzero $R(G)_\mathfrak{p}$-module? The answer can help us identify which orbits in the space contribute to the equivariant $K$-theory after localization at $\mathfrak{p}$, and is given by the following
\begin{theorem}\cite[Proposition 3.7(iv)]{Se}\label{specsupp}
Let $\mathfrak{p}$ be a prime ideal of $R(G)$. Define the support of $\mathfrak{p}$ to be the minimal topologically cyclic subgroup $A$ of $G$, up to conjugation, such that $\mathfrak{p}$ lies in the preimage of a prime ideal of $R(A)$ under the restriction map $R(G)\to R(A)$. Let $H$ be a Lie subgroup of $G$. Then $R(H)_\mathfrak{p}$ as a $R(G)_\mathfrak{p}$-module is not zero if and only if $A$ is conjugate to a subgroup of $H$. 
\end{theorem}
This representation theoretic result is best understood in the following geometric terms. The representation ring $R(G)$ can be thought of as the ring of characters, which are a special kind of class functions on $G$. In parallel to recovering a locally compact Hausdorff space from the ring of its continuous complex-valued functions by taking the set of maximal ideals, we consider $\text{Spec }R(G)$, the algebro-geometric incarnation of $G$.\footnote{The spectrum of the complexified representation ring is the GIT quotient of the complexified group $G_\mathbb{C}$ by the conjugation action by itself, whose points correspond to the semi-simple conjugacy classes of $G_\mathbb{C}$ (\cite[\S 1.1]{FHT})} Note that the affinization functor $G\to \text{Spec }R(G)$ is covariant. As a $R(G)_\mathfrak{p}$-module, $R(H)_\mathfrak{p}$ is not zero if and only if the stalk of the sheaf of $R(G)$-modules $\mathcal{R}(H)$ associated with $R(H)$ at $\mathfrak{p}$ is non-empty,  or equivalently the image of the map $\text{Spec }R(H)\to \text{Spec }R(G)$ induced by the embedding $H\hookrightarrow G$ contains the point $\mathfrak{p}$. By definition the support $A$ of $\mathfrak{p}$ is the minimal topologically cyclic subgroup, up to conjugation, such that the image of the map $\text{Spec }R(A)\to \text{Spec }R(G)$ contains $\mathfrak{p}$. So Theorem \ref{specsupp} says that the image of $\text{Spec }R(H)\to\text{Spec }R(G)$ contains $\mathfrak{p}$ if and only if $H$ contains $A$ up to conjugation.
\begin{theorem}[Segal localization theorem, general version]\label{segallocg}
	Let $\mathfrak{p}$ be a prime ideal of $R(G)$, $A$ a support of $\mathfrak{p}$, and $X$ a locally compact Hausdorff $G$-space. Then the restriction map
	\[i^*\colon K_G^*(X)\to K_G^*(G\cdot X^A)\]
	becomes an isomorphism after localizing the coefficient ring $R(G)$ at $\mathfrak{p}$.
\end{theorem}
Note that $G\cdot X^A$ comprises precisely those orbits whose stabilizers contain conjugates of $A$. By repeating the arguments in the proof of the abelian version and applying Theorem \ref{specsupp}, we see that on localizing at $\mathfrak{p}$ the contribution to equivariant $K$-theory only comes from $G\cdot X^A$, and hence the result.
\begin{remark}
	Let $S$ be a multiplicative subset of $R(G)$ and $i_x^*\colon R(G)\to R(G_x)$ the restriction map to the representation ring of the stabilizer subgroup $G_x$. In analogy with the setting in equivariant cohomology, we define
	\[X^S=\{x\in X|\ i_x^*(\chi)\neq 0\text{ for all }\chi\in S\}\]
	(cf. \cite[Definition 3.1.3]{AP}). Take $S=R(G)\setminus \mathfrak{p}$. Then by Theorem \ref{specsupp}, $X^{R(G)\setminus\mathfrak{p}}=G\cdot X^A$, and Theorem \ref{segallocg} can be rephrased as the restriction map $K_G^*(X)_\mathfrak{p}\to K_G^*(X^{R(G)\setminus\mathfrak{p}})_\mathfrak{p}$ being an isomorphism. This formulation of the Segal localization theorem is analogous to the localization theorem for equivariant cohomology as couched in \cite[Theorem 3.2.6]{AP}.
\end{remark}
Let us end this section by pointing out how the idea of affinization, together with the localization theorem, gives an algebro-geometric construction of equivariant cohomology theories in general (see \cite{Gro, RK}). The $R(G)$-module $K^*_G(X)$ gives rise to the sheaf of modules $\mathcal{K}_G^*(X)$ over $\text{Spec }R(G)$. Conversely, one may retrieve $\mathcal{K}_G^*(X)$ by stitching the stalks at all prime ideals $\mathfrak{p}$ which are specified by the localization theorem, and recover $K_G^*(X)$ by taking global sections. This `delocalization' approach also applies to Borel equivariant cohomology (\cite[Introduction]{RK}), and was used in \cite{Gro} to define equivariant elliptic cohomology. The latter is done by a base change to every stalk of the sheaf of modules associated to equivariant cohomology resulting from comparing the coefficient rings of equivariant cohomology and elliptic cohomology.

\section{Localization formula for equivariant $K$-theory}
This section is devoted to another central result in equivariant $K$-theory, namely the Atiyah-Segal localization formula. The formula is another vivid illustration of the mantra that localization of the coefficient ring helps extract global information from that of the fixed point set. In the setting of the formula the global information concerned is the equivariant index of an elliptic operator, which by the Atiyah-Singer index theorem can be got using pushforward maps (aka Gysin maps and wrong-way maps) for equivariant $K$-theory. As such we spend a major part of this section on setting up the background for pushforward maps before coming to the localization formula in the last subsection. In the rest of this section, we assume that $G$ is a compact Lie group unless otherwise specified.\subsection{Thom isomorphism and Euler classes}
For a real vector bundle $E\to X$ of rank $n>0$, its orientability can be characterized by the existence of a Thom class $\tau_E\in H^n_c(E; \mathbb{Z})$, the cohomology class which restricts consistently to a preferred generator of $H_c^n(\mathbb{R}^n; \mathbb{Z})\cong H^n(S^n; \mathbb{Z})\cong\mathbb{Z}$, which in turn reflects the orientation. This leads to the more general notion of orientability of a vector bundle with respect to a generalized cohomology theory, which is tantamount to the existence of a Thom class in the generalized cohomology in a suitable sense. Let $E$ be a vector bundle of even real rank $2n$. We say $E$ is $K$-orientable if there exists a $K$-theory class $\tau_E\in K^{2n}(E)$, called the Thom class, such that the restriction map $i_x^*\colon K^{2n}(E)\to K^{2n}(E_x)\cong K^{2n}(\mathbb{R}^{2n})$ to each fiber takes $\tau_E$ to $1$ in $K^{2n}(\mathbb{R}^{2n})$, canonically identified with $K^0(\text{pt})\cong\mathbb{Z}$ through the suspension isomorphism. If we identify $K^{2n}(E)$ with $K^0(E)$ through Bott periodicity, then $\tau_E$ restricts to the $K$-theory class in $K^0(\mathbb{R}^{2n})$ which corresponds to $(-\beta)^n\in K^{-2n}(\text{pt})$. If $E$ is of odd rank, then we say it is $K$-orientable if $E\oplus (X\times\mathbb{R})$ is.

It is a well-known result by Atiyah-Bott-Shapiro that $K$-orientability is equivalent to the existence of a $\text{Spin}^c$ structure of the vector bundle, and that a Thom class can be constructed by clutching two spinor bundles through Clifford module multiplication (cf. \cite[Theorem 12.3 and the following remark]{ABS}). This prompts us to define $K$-orientability in the equivariant setting analogously in terms of equivariant $\text{Spin}^c$ structure.\footnote{Defining equivariant $K$-orientability as the existence of an equivariant $K$-theory class of $E$ which restricts to a preferred generator of the equivariant $K$-theory of each fiber of $E$, in analogy to the definition of ordinary $K$-orientability, does not make sense, because the restriction map is not equivariant if a fiber is not invariant under the $G$-action. So one has to resort to equivariant $\text{Spin}^c$-structures to define equivariant $K$-orientability.} 
\begin{definition}\label{eqkor}
	Let $E$ be an oriented real $G$-vector bundle of rank $n$ over $X$ where the $G$-action preserves the orientation. Then $E$ has an equivariant $\text{Spin}^c$ structure if its oriented orthonormal frame bundle $P$ (with respect to a fiberwise inner product on $E$ invariant under the $G$-action, which always exists by the averaging trick) can be lifted to a principal $\text{Spin}^c(n)$-bundle $\widetilde{P}$ and the $G$-action on $P$ can be lifted to one on $\widetilde{P}$. If $E$ admits an equivariant $\text{Spin}^c$-structure then it is said to be equivariantly $K$-orientable. If $X$ is compact, then the equivariant Thom class $\tau_E$ corresponding to an equivariant $\text{Spin}^c$ structure is defined to be the ordinary Thom class of $E$ constructed using the spinor bundle associated with the $\text{Spin}^c$ structure and equipped with the $G$-action coming from the lifted $G$-action on $\widetilde{P}$.
\end{definition}
\begin{remark}
\begin{enumerate}
\item An equivariant Thom class depends on the equivariant $\text{Spin}^c$ structure which is used to construct the spinor bundles needed to define the Thom class itself. The set of equivariant $\text{Spin}^c$ structures is in bijective correspondence with that of integral lifts of the second Stiefel-Whitney class of the vector bundle, and is acted upon by the Picard group of equivariant complex line bundles over $X$: the equivariant principal $\text{Spin}^c$-bundle $\widetilde{P}$ is taken by the equivariant line bundle $L$ to another equivariant principal $\text{Spin}^c$-bundle obtained by scaling the transition functions of $\widetilde{P}$ by those of the circle bundle of $L$.
\item Any equivariant complex vector bundle admits a natural equivariant $\text{Spin}^c$-structure coming from the hermitian structure through the homomorphism of structure groups $U(n)\to\text{Spin}^c(2n)$ (see \cite[end of \S 3]{ABS}), and hence is equivariantly $K$-orientable. 
\item Atiyah showed in \cite{At2} using index theoretic arguments that there is also the Thom isomorphism in the equivariant setting, realized by the multiplication by an equivariant Thom class. This justifies the definition of equivariant $K$-orientability as the existence of equivariant $\text{Spin}^c$-structures.
\end{enumerate}
\end{remark}
\begin{theorem}[Equivariant Thom isomorphism]
	Let $E$ be an equivariant $\text{Spin}^c$ $G$-vector bundle over $X$. Then the map
	\begin{align*}
		i_*\colon K_G^*(X)&\to K_G^{*+\text{rank }E}(E)\\
		V&\mapsto \pi^*V\cdot\tau_E
	\end{align*}
	is an isomorphism of $K_G^*(X)$-modules.
\end{theorem}
The Euler class for equivariant $K$-theory can be defined in analogy with the cohomological Euler class as follows. The equivariant $K$-theoretic Euler class also measures the `twistedness' of equivariant vector bundles like its cohomological counterpart does, and make an appearance in the localization formula as we will see.\footnote{The `twistedness' here measures two aspects: how far the vector bundle, without regard to the equivariant structure, is from being trivial, as well as how nontrivial the group action is.}
	\begin{definition}\label{euler}
		Let $E\to X$ be a $K$-oriented $G$-vector bundle, $i^*\colon K_G^*(E)\to K_G^*(X)$ the restriction of $E$ to its zero section. The equivariant Euler class $e(E)$ is defined to be $i^*\tau_E$.
	\end{definition}
	The Thom class is multiplicative in the sense that for two equivariantly $K$-oriented vector bundles $E_1$ and $E_2$, we have $\pi_1^*\tau_{E_1}\cdot\pi_2^*\tau_{E_2}=\tau_{E_1\oplus E_2}$, where $\pi_i$ is the projection from $E_1\oplus E_2$ to $E_i$. By Definition \ref{euler} the Euler class is multplicative as well.
	\begin{example}\label{CP}
		Consider the complex line bundle $\mathbb{C}_n\to\text{pt}$, on which $S^1$ acts with weight $n$. We would like to work out $\tau_{\mathbb{C}_n}$ and $e(\mathbb{C}_n)$ (associated to the equivariant $\text{Spin}^c$ structure coming from the complex structure). We shall first compute $K_{S^1}^*(\mathbb{C}_n^+)=K_{S^1}^*(\mathbb{CP}^1)$, where we identify $\mathbb{C}_n^+$ with $\mathbb{CP}^1$ through $z\mapsto [1:z]$, $\infty\mapsto[0:1]$, and $S^1$ acts by $e^{i\theta}\cdot[z_0:z_1]=[z_0:e^{in\theta}z_1]$. Consider the closed cover $D_{[1:0]}\cup D_{[0:1]}$, where $D_{[1:0]}:=\{[1: z]| |z|\leq 1\}$ and $D_{[0:1]}:=\{[z:1]| |z|\leq 1\}$. The intersection $D_{[1:0]}\cap D_{[0:1]}$ is the circle $\mathcal{S}:=\{[1:z]| |z|=1\}$. We have the following Mayer-Vietoris sequence with respect to the closed cover.
		\begin{eqnarray*}
			\xymatrix@R-0.5pc{0\ar[r]&K_{S^1}^0(\mathbb{CP}^1)\ar[r]& K_{S^1}^0(D_{[1:0]})\oplus K_{S^1}^0(D_{[0:1]})\ar[r]\ar[d]^\cong&K_{S^1}^0(\mathcal{S})\ar[r]\ar[d]^\cong&\cdots\\ &&R(S^1)\oplus R(S^1)\ar[d]^\cong\ar[r]& R(\mathbb{Z}_n)\ar[d]^\cong&\\ &&\mathbb{Z}[t, t^{-1}]\oplus\mathbb{Z}[t, t^{-1}]\ar[r]^{r}&\mathbb{Z}[u]/{(u^n-1)}&}
		\end{eqnarray*}
		This sequence merits some explanations. As both $D_{[0:1]}$ and $D_{[0:1]}$ are equivariantly contractible, their equivariant $K$-theory rings are isomorphic to $K_{S^1}^*(\text{pt})\cong R(S^1)$, resulting in the left vertical isomorphism. The element $t$ is represented by the trivial line bundles over $D_{[1:0]}$ and $D_{[0:1]}$ with $S^1$ acting on the fibers with weight 1. On the other hand, $\mathcal{S}$ is acted upon by $S^1$ transitively with the stabilizer subgroup $\mathbb{Z}_n$. By Example \ref{ekex} (\ref{ekhomog}), $K_{S^1}^*(\mathcal{S})=K_{S^1}^*(S^1/\mathbb{Z}_n)\cong K_{\mathbb{Z}_n}^*(\text{pt})\cong R(\mathbb{Z}_n)\cong \mathbb{Z}[u]/(u^n-1)$, where $u$ is represented by the trivial line bundle over $\mathcal{S}$ with $S^1$ acting on the fibers with weight 1. This gives the right vertical isomorphisms. The left end of the sequence, which is supposedly $K_{S^1}^{-1}(\mathcal{S})$, is then 0. The restrictions of $K_{S^1}^*(D_{[1:0]})$ and $K_{S^1}^*(D_{[0:1]})$ to $K_{S^1}^*(\mathcal{S})$ send $t$ to $u$. Thus the bottom map $r$, identified through the vertical isomorphisms with the signed sum of restrictions from $D_{[1:0]}$ and $D_{[0:1]}$ to $\mathcal{S}$, sends a pair of Laurent polynomials $(p(t), q(t))$ to $p(u)-q(u)$ modulo $u^n-1$. By exactness of the sequence, the surjectivity of $r$ and the vanishing of $K_{S^1}^{1}(\mathcal{S})$, $K_{S^1}^1(D_{[1:0]})$ and $K_{S^1}^1(D_{[0:1]})$, we have $K^{1}_{S^1}(\mathbb{CP}^1)=0$. So $K_{S^1}^*(\mathbb{CP}^1)$ is isomorphic to $\text{ker}(r)$, which is the $R(S^1)$-algebra $\{(p(t), q(t))\in\mathbb{Z}[t, t^{-1}]\oplus\mathbb{Z}[t, t^{-1}]|t^n-1\text{ divides }p(t)-q(t)\}$. In fact, it is a free $R(S^1)$-algebra generated by $(1-t^n, 0)$ and $(1, 1)$, and thus its rank is 2, verifying the claim we obtain immediately after Theorem \ref{segallocab}. It follows that \begin{align*}K_{S^1}^*(\mathbb{C}_n)=\widetilde{K}^*_{S^1}(\mathbb{CP}^1)&\cong\{(p(t), 0)\in\mathbb{Z}[t, t^{-1}]\oplus\mathbb{Z}[t, t^{-1}]| t^n-1\text{ divides }p(t)\}\\ &\cong R(S^1)\cdot (1-t^n, 0).\end{align*}
	It can be checked that for the ordinary Thom class of $\mathbb{C}$, which by definition is $1-H\in \widetilde{K}^0(\mathbb{CP}^1)$, its equivariant lifts in $\widetilde{K}^*_{S^1}(\mathbb{CP}^1)$ must be $(t^m-t^{m+n}, 0)$ for $m\in\mathbb{Z}$ ($S^1$ acts on the fibers of the trivial line bundle with weight $m$, $H|_{[1:0]}$ with weight $m+n$ and $H|_{[0:1]}$ with weight $m$). All these lifts, differing from each other by the module multiplication by some equivariant line bundle over a point, are possible equivariant Thom classes. However, there is only one lift among them which corresponds to the equivariant $\text{Spin}^c$ structure stemming from the complex structure of $\mathbb{C}_n$, namely the generator $(1-t^n, 0)$. We have that the equivariant Thom class $\tau_{\mathbb{C}_n}$ is $(1-t^n, 0)$ by Definition \ref{eqkor}. It follows that the equivariant Euler class is $e(\mathbb{C}_n)=1-t^n\in K_{S^1}^*(\text{pt})$. 
	
	More generally, one can show similarly that if a compact torus $T$ acts on $\mathbb{C}^m_{\alpha_1, \cdots, \alpha_m}$, which is decomposed into weight spaces $\bigoplus_{i=1}^m\mathbb{C}_{\alpha_i}$, $\alpha_i\in\Lambda^*_T$, then by multiplicativity $ e(\mathbb{C}^m_{\alpha_1, \cdots, \alpha_m})=\prod_{i=1}^m(1-t^{\alpha_i})\in K_T^*(\text{pt})$. This equivariant Euler class is not zero in the coefficient ring if and only if all the isotypical weights are nonzero, despite the vanishing of the ordinary Euler class due to the triviality of $\mathbb{C}^m_{\alpha_1, \cdots, \alpha_m}$ as an ordinary vector bundle. This distinction between equivariant and ordinary Euler classes is an important feature which sets equivariant $K$-theory (and other generalized equivariant cohomology theories) apart from its ordinary counterpart, and forms the backbone of the proof of the localization formula.
	
	If we expand $ \prod_{i=1}^m(1-t^{\alpha_i})$ and write as a formal sum of vector bundles, we have $ e(\mathbb{C}^m_{\alpha_1, \cdots, \alpha_m})=\sum_{i=0}^m (-1)^i\bigwedge\nolimits^i\mathbb{C}^m_{\alpha_1, \cdots, \alpha_m}$. In fact, for a general equivariant complex vector bundle $E$, its equivariant Euler class is the alternating sum of its exterior powers $\sum_{i=0}^{\text{rank }E}(-1)^i\bigwedge\nolimits^i E$ (\cite[p. 140, 3rd paragraph]{Se2}).
	\end{example}
\subsection{The pushforward maps}
\begin{definition}\label{mapkor}
	Let $f\colon X\to Y$ be an equivariant map between $G$-manifolds. The map $f$ is said to be equivariantly $K$-orientable if there exists an equivariant $\text{Spin}^c$ $G$-vector bundle $E$ over $Y$ and an equivariant embedding $i\colon X\to E$ such that the normal bundle $N_{E/X}$ of $X$ in $E$ is equivariantly $\text{Spin}^c$ and $f$ factors as $\pi\circ i$, where $\pi$ is the projection $E\to Y$.
\end{definition}
A necessary condition for $f$ to be equivariantly $K$-orientable is that $TX\oplus f^*TY$ be equivariantly $\text{Spin}^c$. This condition is also sufficient for $f$ to be equivariantly $K$-orientable if in addition $X$ has finite orbit type (e.g. $X$ is compact). That is because in this case $X$ can be equivariantly embedded into a real linear representation $V$ of $G$ (cf. \cite{Mo} and \cite{Pa}), $E$ then can be taken to be $Y\times (V\oplus V)$ which can be made a complex vector bundle in a natural way and thus an equivariantly $\text{Spin}^c$ vector bundle over $Y$, and the normal bundle of $X$ for the embedding $i\colon X\to Y\times (V\oplus V)$ can be seen to be equivariantly $\text{Spin}^c$.
\begin{definition}\label{pushforward}
	\begin{enumerate}
		\item\label{embedding} Let $i\colon X\to Y$ be an equivariant embedding such that the normal bundle $N_{Y/X}$ is equivariantly $\text{Spin}^c$. The pushforward map $i_*\colon K_G^*(X)\to K_G^{*+\text{dim }Y-\text{dim }X}(Y)$ is defined to be the composition of maps
		\[K_G^*(X)\to K_G^{*+\text{dim }Y-\text{dim }X}(N_{Y/X})\to K_G^*(U)\to K_G^*(Y),\]
		where the first map is the Thom isomorphism, the second induced by the equivariant diffeomorphism between $N_{Y/X}$ and an open tubular neighborhood $U$ of $X$, the third induced by the map $Y\to U^+$ which collapses $Y\setminus U$ to the point of compactification.
		\item The pushforward map $\pi_*\colon K_G^*(E)\to K_G^{*-\text{rank }E}(Y)$ induced by the projection from an equivariantly $\text{Spin}^c$ vector bundle $E$ to $Y$ is defined to be the inverse of the Thom isomorphism $i_*\colon K_G^*(Y)\to K_G^{*+\text{rank }E}(E)$.
		\item If $f\colon X\to Y$ is equivariantly $K$-orientable with $E$, $i$ and $\pi$ as in Definition \ref{mapkor}, then the induced pushforward map $f_*\colon K_G^*(X)\to K_G^{*+\text{dim }Y-\text{dim }X}(Y)$ is defined to be $\pi_*\circ i_*$. 
	\end{enumerate}
\end{definition}
The pushforward map $f_*$ in fact only depends on the equivariant $\text{Spin}^c$ structure of $TX\oplus f^*TY$ but not on the choice of $E$. 
\begin{example}
	\begin{enumerate}
		\item The Thom isomorphism $i_*\colon K_G^*(X)\to K_G^*(E)$ is the pushforward map induced by the inclusion of $X$ into $E$ as the zero section. 
		\item Let $X$ be a compact and equivariantly $\text{Spin}^c$ manifold of dimension $2n$. Then the collapsing map $\pi_X\colon X\to \text{pt}$ is equivariantly $K$-orientable. The pushforward $\pi_{X*}(E)\in K_G^{-2n}(\text{pt})\cong R(G)$ is called the equivariant topological index of $E$. On the other hand, one can define the Dirac operator 
		\[\slashed{\partial}\colon \Gamma(S^+\otimes E)\to \Gamma(S^-\otimes E)\]
		where $S^+$ and $S^-$ are the positive and negative spinor bundles. The equivariant analytic index, which is $\text{ker}(\slashed{\partial})-\text{coker}(\slashed{\partial})$, is a finite dimensional virtual $G$-representation. The Atiyah-Singer index theorem asserts the equality of the equivariant topological and analytic indices (cf. \cite{AS}), and explains the significance of $\text{Spin}^c$ structures in $K$-theory from the analytic perspective. If $M$ is further assumed to be a complex $G$-manifold with the equivariant $\text{Spin}^c$ structure coming from the complex structure, then the Dirac operator is just (a certain constant multiple of) the Dolbeault operator $\overline{\partial}$ on the rolled-up Dolbeault complex twisted by $E$. The equivariant analytic index then is the holomorphic Euler characteristic $\chi(M, E):=\sum_i(-1)^i H^i(M, E)$. 
		\item Let $X$ be a compact Hamiltonian $G$-manifold with a symplectic form $\omega$ which defines an integral cohomology class, and $L\to X$ a prequantum line bundle, which satisfies $c_1(L)=[\omega]$. With the moment map, one can lift the $G$-action on $X$ to one on $L$ by Kostant's formula. The geometric quantization of $L$ is defined to be its equivariant analytic index, which by the index theorem is $\pi_{X*}(L)$ with respect to the $\text{Spin}^c$-structure coming from the almost complex structure of $X$ (\cite{Me}). This generalizes the geometric quantization of a K\"ahler $G$-manifold defined using the complex polarization of a holomorphic prequantum line bundle. Let $X_0$ be the symplectic reduction at $0\in\mathfrak{g}^*$, and $L_0$ be the prequantum line bundle over $X_0$ which is the reduction of $L$. The principle of `quantization commutes with reduction' says that $\pi_{X_0*}(L_0)$ is the dimension of the trivial subrepresentation of $\pi_{X*}(L)$ (\cite[Theorem 1.1]{Me}).% In other words, the quantum analogue of symplectic reduction is taking the invariant subspace of the group action.
	\end{enumerate}
\end{example}
\subsection{Atiyah-Segal localization formula}
We are now in a position to state the Atiyah-Segal localization formula, which expresses the equivariant topological index as a sum of indices of fixed point sets. For simplicity we state the formula for equivariant complex manifolds.
	\begin{theorem}[Atiyah-Segal localization formula for complex manifolds \cite{AtSe}]
		Let $T$ be a compact torus, $X$ a compact complex $T$-manifold and $E$ a holomorphic $T$-vector bundle. Then
		\[\pi_{X*}(E)=\sum_{F\subseteq X^T}\pi_{F*}\left(\frac{i_F^*E}{e(N^*_{X/F})}\right),\]
		where $F$ ranges over the connected components of $X^T$, $i_F: F\hookrightarrow X$ the inclusion of $F$ into $X$, and $\pi_F$ and $\pi_X$ the collapsing maps on $F$ and $X$ respectively. 
		%\begin{eqnarray*}
		%\xymatrix{K_T^*(F)\ar[r]^{i_{F*}}\ar[dr]_{\pi_{F*}}& K_T^*(X)\ar[d]^{\pi_{X*}}\ar[r]^{i_F^*}& K_T^*(F)\\ & K_T^*(\text{pt})&}
	%\end{eqnarray*}
		\end{theorem}
	\begin{eqnarray*}
		\xymatrix@+3.5pc{\bigoplus_{F\subseteq X^T}K_T^*(F)\ar[r]^{\bigoplus_{F\subseteq X^T}i_{F*}}\ar[dr]_{\bigoplus_{F\subseteq X^T}\pi_{F*}}& K_T^*(X)\ar[d]^{\pi_{X*}}\ar[r]^{\bigoplus_{F\subseteq X^T}i_F^*}& \bigoplus_{F\subseteq X^T}K_T^*(F)\\ & K_T^*(\text{pt})&}
	\end{eqnarray*}
		It is known that $F$ as a $T$-fixed submanifold is a complex submanifold, and that the normal bundle $N_{X/F}$ is an equivariant complex vector bundle. If $\pi_{F*}$ and $\pi_{X*}$ are $K$-oriented by the (inherited) complex structure, then the $K$-orientation of $i_{F*}$ compatible with the equation $\pi_{F*}=\pi_{X*}\circ i_{F*}$ is induced by the \emph{opposite} of the inherited complex structure of $N_{X/F}$. 	%\begin{eqnarray*}
%		\xymatrix@R-1pc{K_T^*(F)\ar[r]^{i_{F*}}\ar[dr]_{\pi_{F*}}& K_T^*(X)\ar[d]^{\pi_{X*}}\ar[r]^{i_F^*}& K_T^*(F)\\ & K_T^*(\text{pt})&}
%	\end{eqnarray*}
	By Definitions \ref{euler} and \ref{pushforward} (\ref{embedding}), $i_F^*\circ i_{F*}$ is the multiplication by $e(N^*_{X/F})$. On each fiber of $N_{X/F}$, there is no zero isotypical weight, for otherwise the zero weight subspace would also be fixed by $T$ and hence part of $F$. The same is true of $N^*_{X/F}$. This leads to the distinctive feature of $e(N^*_{X/F})$ that it is invertible in the localized equivariant $K$-theory $K_T^*(F)_{(0)}\cong R(T)_{(0)}\otimes K^*(F)$, which in general is not true in the nonequivariant setting. 
	
	Now more explanations are in order for this feature of the equivariant Euler class. Let $ N^*_{X/F}=\sum_i t^{\alpha_i}\otimes E_{\alpha_i}\in R(T)\otimes K^*(F)$ be the isotypical decomposition as in Example \ref{ekex} (\ref{isotypical}). Then by multiplicativity and the remark at the end of Example~\ref{CP}, 
	\begin{align*}
		&e(N^*_{X/F})=\prod_{i}e(t^{\alpha_i}\otimes E_{\alpha_i})=\prod_i\left(\sum_{j=0}^{\text{rank }E_{\alpha_i}}(-1)^jt^{j\alpha_i}\otimes\bigwedge\nolimits^j E_{\alpha_i}\right)\\
					=&\prod_i\left((1-t^{\alpha_i})^{\text{rank }E_{\alpha_i}}\otimes 1+\sum_{j=0}^{\text{rank }E_{\alpha_i}}(-1)^jt^{j\alpha_i}\otimes\left(\bigwedge\nolimits^j E_{\alpha_i}-\binom{\text{rank }E_{\alpha_i}}{j}\right)\right).
	\end{align*}
	The $K$-theory class $\bigwedge\nolimits^j E_{\alpha_i}-\binom{\text{rank }E_{\alpha_i}}{j}$ lives in the reduced $K$-theory of the compact space $F$, and so is nilpotent in $K^*(F)$ (\cite[3.1.6]{At}).\footnote{The compact manifold $F$ can be covered by finitely many contractible closed sets $\{A_i\}_{i=1}^n$. The reduced $K$-theory $\widetilde{K}^*(F)$ is naturally isomorphic to ${K}^*(F, A_i)$ for each $i$. So given any $n$ $K$-theory classes $x_1, \cdots, x_n\in \widetilde{K}^*(F)$, each of them can be identified with a relative $K$-theory class of ${K}^*(F, A_i)$, and their product then lives in ${K}^*(F, \bigcup_{i=1}^n A_i)={K}^*(F, F)$, which is 0.} The term $\sum_{j=0}^{\text{rank }E_{\alpha_i}}(-1)^j t^{j\alpha_i}\otimes\left(\bigwedge\nolimits^j E_{\alpha_i}-\binom{\text{rank }E_{\alpha_i}}{j}\right)$ then is nilpotent in $R(T)\otimes K^*(F)$ as well. Now each factor is a sum of a nilpotent element and $(1-t^{\alpha_i})^{\text{rank }E_{\alpha_i}}\otimes 1$, which is invertible in $R(T)_{(0)}\otimes K^*(F)$ as $\alpha_i$ is nonzero. It follows that each factor and their product are invertible. 
	
	After localizing the first row of the above diagram at the zero ideal, the map $\bigoplus_{F\subseteq X}i_F^*$ becomes an isomorphism by the Segal localization theorem. So is the composition $\left(\bigoplus_{F\subseteq X}i_{F*}\right)\circ\left(\bigoplus_{F\subseteq X} i^*_F\right)$ as it is the componentwise multiplication by $e(N_{X/F}^*)$, which is invertible. The pushforward $\bigoplus_{F\subseteq X}i_{F*}$ then is an isomorphism as well. So one can rewrite $\pi_{X*}(E)$ as the composition $\pi_{X*}\circ \left(\bigoplus_{F\subseteq X}i_{F*}\right)\circ\left(\bigoplus_{F\subseteq X}i_F^*\circ \bigoplus_{F\subseteq X}i_{F*}\right)^{-1}\circ \left(\bigoplus_{F\subseteq X}i_F^*E\right)$, culminating in the desired localization formula.
	\begin{remark}
		There is also the same localization formula for $T$-$\text{Spin}^c$ manifolds of even dimension, but one has to equivariantly $K$-orient the pushforward $\pi_{X*}$, $\pi_{F*}$ and the normal bundle $N_{X/F}$ involved in the formula with suitable equivariant $\text{Spin}^c$ structures so that $\pi_{F*}=\pi_{X*}\circ i_{F*}$ holds. For instance, the normal bundle, being equipped with a torus action which fixes only the zero section, can be made an equivariant complex vector bundle and $K$-oriented by the associated equivariant $\text{Spin}^c$ structure (\cite[Lemma 2.1]{HL}). Together with the given equivariant $\text{Spin}^c$ structure of the ambient manifold, one can $K$-orient the fixed point set in a compatible manner.
	\end{remark}
	To evaluate the ordinary index of a manifold, very often its symmetry can be exploited to our advantage through the application of the localization formula, as it reduces to the computation of indices of the fixed point set which is easier to deal with. This is illustrated by the next example.
	\begin{example}\label{topindH}
		The ordinary topological index of the hyperplane bundle $H$ of $\mathbb{CP}^1$ can be got by computing its equivariant topological index with respect to, for instance, the $S^1$-action on $\mathbb{CP}^1$ and $H$ as in Example \ref{CP}, and applying the forgetful map which amounts to replacing each irreducible character to its dimension, namely 1. Recall that the weights of the $S^1$-action on the normal bundles of the two fixed points $[1:0]$ and $[0:1]$, and the fibers $H|_{[1:0]}$ and $H|_{[0:1]}$ are respectively $n$, $-n$, $m+n$ and $m$ for $m, n\in\mathbb{Z}$. The localization formula gives the equivariant topological index as 
		\[\frac{t^{m+n}}{1-t^{-n}}+\frac{t^m}{1-t^n}=t^m(1+t^n).\]
		Thus the ordinary topological index is 2. This is consistent with the analytic index $\text{dim }H^0(\mathbb{CP}^1, H)=2$. The generator of $\widetilde{K}^*(\mathbb{CP}^1)$ with ordinary topological index 1 then is $H-1=1-H^{-1}$, which is taken to be the Bott class.
	\end{example}
\begin{example}
		The localization formula offers an interesting proof of the Weyl character formula for representations of compact Lie groups. Let $G$ be a connected compact Lie group and $T$ its maximal torus, which acts on the generalized flag manifold $G/T$ by left translation. Fix a complex structure of $G/T$ by specifying the set of positive roots $R^
		+$. Consider $\pi_*\colon K_T^*(G/T)\to K_T^*(\text{pt})$ induced by this complex structure and $E_\mu:=G\times_T\mathbb{C}_\mu$ where $\mu$ is a dominant weight. The equivariant analytic index $\chi(G/T, E_\mu)=H^0(G/T, E_\mu)$ is the irreducible representation of $G$ with highest weight $\mu$ by the Borel-Weil-Bott Theorem, while the equivariant topological index $\pi_*(E_\mu)$ can be computed using the Atiyah-Segal localization formula as follows. The fixed point set $(G/T)^T$ is $\{wT/T|w\in W\}$ where $W$ is the Weyl group. The normal bundle of the identity coset decomposes into positive root spaces, and so $N^*_{(G/T)/(wT/T)}\cong\bigoplus_{\alpha\in R^+}\mathbb{C}_{-w\alpha}$, and $e(N^*_{(G/T)/(wT/T)})=\prod_{\alpha\in R^+}(1-t^{-w\alpha})$. The restriction $i_{wT/T}^*E_\mu$ is $t^{w\mu}$. Assembling these data into the localization formula then yields
		\[\pi_{*}(E_\mu)=\sum_{w\in W}\frac{t^{w\mu}}{\prod_{\alpha\in R^+}(1-t^{-w\alpha})}, \]
		which is the Weyl character formula for the irreducible representation with highest weight $\mu$.
	\end{example}

\section{Applications to equivariant formality}
Throughout this section, we assume that $G$ is a compact Lie group unless otherwise specified.
\subsection{Equivariant formality in $K$-theory} Equivariant formality, first identified and discussed at length in \cite{BBFMR} and named in \cite{GKM}, is an important property of group actions on topological spaces which allows for easy computation of their (equivariant) cohomology. A $G$-action on a space $X$ is said to be equivariantly formal if the Leray-Serre spectral sequence for the rational cohomology of the fiber bundle $X\hookrightarrow X\times_G EG\to BG$ collapses on the $E_2$-page. Examples of interest which are known to be equivariantly formal abounds, e.g. Hamiltonian group actions on compact symplectic manifolds, linear algebraic torus actions on smooth complex projective varieties, any space with a CW decomposition invariant under the group action (\cite[\S 1.2 and Theorem 14.1]{GKM}), and conjugation action on any compact Lie group by itself (\cite[\S 11.9, Item 6]{GS}, \cite{J}). Let $X$ be a $G$-space where $G$ is a compact Lie group. The desirability of equivariant formality as a property of group actions can be seen from the following characterizations of equivariant formality of the $G$-action on $X$.
\begin{enumerate}
	\item The equivariant cohomology $H_G^*(X; \mathbb{Q})$ is isomorphic to 
	\[H_G^*(\text{pt}; \mathbb{Q})\otimes H^*(X; \mathbb{Q})\]
	as $H_G^*(\text{pt}; \mathbb{Q})$-modules.  
	\item The forgetful map $H_G^*(X; \mathbb{Q})\to H^*(X; \mathbb{Q})$ is surjective.
	\item\label{dimcount} If $G$ is a torus $T$, then the above conditions are all equivalent to \[\text{dim }H^*(X; \mathbb{Q})=\text{dim }H^*(X^T; \mathbb{Q}).\]
\end{enumerate}
The above three conditions and equivariant formality are equivalent if $X$ is sufficiently nice, e.g. $X$ is a finite dimensional $G$-CW complex with finitely many orbit types (\cite[Theorem 3.10.4]{AP}). The first condition implies that the equivariant cohomology can be recovered from the ordinary cohomology as the former is a free module over the coefficient ring with module basis corresponding to the basis of the latter. The third condition says that equivariantly formal torus actions have the maximal possible fixed point set in the cohomological sense, as in general we have $\text{dim }H^*(X^T; \mathbb{Q})\leq\text{dim }H^*(X; \mathbb{Q})$ (\cite[IV 5.5, p.62]{BBFMR}). %By applying the forgetful map which maps the Leray-Serre spectral sequence of the fiber bundle $X\hookrightarrow X\times_G EG\to BG$ to that of $X$ regarded as a fiber bundle over a single point, we see that the second condition follows from the definition of equivariant formality. also implies it and hence the first condition by the Leray-Hirsch Theorem. 
The second condition asserts that for equivariantly formal actions, any rational ordinary cohomology class admits an equivariant lift. When $X$ is a compact $G$-manifold, this property makes it amenable to the application of the Atiyah-Bott localization formula (\cite{AB2}) when computing the integral of any top degree ordinary differential form. As noted in the Introduction, the Atiyah-Bott localization formula in fact is inspired by the Atiyah-Segal localization formula. So it is tempting to put forth an analogue of equivariant formality for equivariant $K$-theory which makes it possible to compute the ordinary topological index of any vector bundle by the Atiyah-Segal localization formula. In fact the notion of equivariant formality in $K$-theory was first introduced and explored in \cite{HL}, where they coined the term `weak equivariant formality' and exploited the notion to show equivariant formality of Hamiltonian group actions on compact symplectic manifolds.\footnote{According to \cite{HL}, the term `weak' is in reference to this condition (recalled in Definition \ref{weakequivform}) being weaker than $K_G^*(X)\cong R(G)\otimes K^*(X)$, the $K$-theoretic analogue of Condition (1) above. However, weak equivariant formality is a stronger condition than its rationalized version that $K_G^*(X)\otimes_{R(G)}\mathbb{Q}\to K^*(X)\otimes\mathbb{Q}$ is an isomorphism, which is equivalent to equivariant formality by Theorem \ref{RKEF}. So this terminology is a misnomer in retrospect.} 
\begin{definition}\label{weakequivform}
	Let $G$ be a compact Lie group and $X$ a locally compact Hausdorff $G$-space. We say $X$ is a weakly equivariantly formal $G$-space if the forgetful map $f_G\colon K_G^*(X)\to K^*(X)$ induces an isomorphism $K_G^*(X)\otimes_{R(G)}\mathbb{Z}\to K^*(X)$, where the $R(G)$-module structure of the tensor factor $\mathbb{Z}$ is induced by the augmentation map $R(G)\to\mathbb{Z}$.
\end{definition}
Inspired by Condition (2) and weak equivariant formality, in \cite{F} we defined the related notion of \emph{rational K-theoretic equivariant formality} as follows.
\begin{definition}
	Let $G$ be a compact Lie group and $X$ a locally compact Hausdorff $G$-space. We say $X$ is a rational $K$-theoretic equivariantly formal (RKEF for short) $G$-space if the forgetful map over the rationals
	\[f_G\otimes\text{Id}_\mathbb{Q}\colon K_G^*(X)\otimes\mathbb{Q}\to K^*(X)\otimes\mathbb{Q}\]
	is onto.
\end{definition}
It turns out that both RKEF and weak equivariant formality over the rationals are the right notions of equivariant formality in the $K$-theoretic sense.
\begin{theorem}[\cite{F}]\label{RKEF} Let $G$ be a compact connected Lie group which acts on a finite CW-complex $X$. Then the following are equivalent.
\begin{enumerate}
	\item $X$ is a RKEF $G$-space. 
	\item $X$ is an equivariantly formal $G$-space. 
	\item $X$ is a weakly equivariantly formal $G$-space over the rationals.
\end{enumerate}
\end{theorem}
Note that, as hinted at in Definition \ref{weakequivform}, the forgetful map $f_G\otimes \text{Id}_\mathbb{Q}\colon K_G^*(X)\otimes\mathbb{Q}\to K^*(X)\otimes \mathbb{Q}$ factors through $K_G^*(X)\otimes_{R(G)}\mathbb{Q}$ and induces the map $K_G^*(X)\otimes_{R(G)}\mathbb{Q}\to K^*(X)\otimes\mathbb{Q}$. The equivalence of conditions (1) and (3) in the above theorem means that the induced map is surjective if and only if it is an isomorphism, which is analogous to the case of equivariant cohomology (the equivalence of Conditions (1) and (2) at the beginning of this section).

Among the aforementioned conditions, the RKEF condition is arguably the most natural characterization of equivariant formality as it has a concrete topological interpretation in terms of vector bundles: more precisely, $X$ is a RKEF $G$-space if for every vector bundle $V$ over $X$ or its suspension $(X\times\mathbb{R})^+$, there are natural numbers $p$ and $q$ such that $V^{\oplus p}\oplus\underline{\mathbb{C}}^q$ admits an equivariant $G$-structure. On the other hand, verifying other conditions, e.g. the first definition of equivariant formality and surjectivity of the forgetul map for equivariant cohomology, involves such algebraic gadget and manipulation as spectral sequences and determining the existence of equivariant extensions of ordinary cohomology classes, which can be tedious (if an equivariant extension of a cohomology class is easy to obtain, it is often the case that the cohomology class is already the characteristic class of an equivariant vector bundle, and in this situation it is more natural to use RKEF to verify equivariant formality). 

We shall pause to remark that Theorem \ref{RKEF} is a topological generalization of a result in algebraic geometry (\cite[Theorem 1.1]{Gr}): it is also an assertion of surjectivity, but of the forgetful map from the rational Grothendieck group of $G$-equivariant coherent sheaves on a $G$-scheme $X$ to the corresponding Grothendieck group for ordinary coherent sheaves, where $G$ is a connected reductive algebraic group. Theorem \ref{RKEF} vindicates the expectation (\cite[Introduction]{Gr}) that the $K$-theoretic forgetful map is onto for equivariantly formal topological spaces. In the following examples we will see how Theorem \ref{RKEF} offers a novel way to show equivariant formality of some spaces which were otherwise proved by cohomological means.
\begin{example}[\cite{F}]\label{fokex}
	\begin{enumerate}
		\item\label{oddgen} Consider the conjugation action of a compact connected Lie group $G$ on itself. Over the suspension of $G$ and for a unitary representation $V$ of $G$, one can construct a vector bundle by gluing two copies of the trivial vector bundle $\textbf{V}$ over two cones of $G$ along $G$ through the clutching function specified by $V$ as the group homomorphism $\rho\colon G\to U(n)$. These vector bundles represent the $K$-theory classes $\delta(\rho)\in K^{-1}(G)$, which by \cite[\S II, Theorem 2.1]{Ho} generate the ring $K^*(G)\otimes\mathbb{Q}$, and can be endowed with equivariant structures: explicitly, we can equip $\textbf{V}:=G\times V$ with the $G$-action $h\cdot(g, v)=(hgh^{-1}, \rho(h)v)$, with respect to which the clutching function $\rho$ is equivariant. It follows that $G$ is a RKEF $G$-space and hence an equivariantly formal $G$-space by Theorem \ref{RKEF}.
		\item\label{eqrank} Let $K$ be a connected Lie subgroup of $G$ and $\text{rank }G=\text{rank }K$. The left translation action on the generalized flag manifold $G/K$ by $G$ is well-known to be equivariantly formal, which can be shown by noting that $G/K$ satisfies the sufficient condition for equivariant formality that its odd cohomology vanishes (\cite[Chapter XI, Theorem VII]{GHV}). Alternatively, we have, by the rationalized version of \cite[Theorem 4.2]{Sn} and the remark following it, that the ring $K^*(G/K)\otimes\mathbb{Q}$ is generated by degree 0 $K$-theory classes represented by homogeneous vector bundles $G\times_KV$, where $V$ is a complex representation of $K$. All these vector bundles obviously have the $G$-action on the first factor by left multplication. By Theorem \ref{RKEF}, we again obtain equivariant formality of the left translation on $G/K$ by $G$.
		\item Let $X$ be a finite $G$-CW complex with maximal rank connected isotropy subgroups. Equivariant formality of the $G$-action follows from \cite[Corollary 3.5]{GR}. Alternatively, to get equivariant formality one can invoke Theorem \ref{RKEF} and show that a stable $G$-equivariant structure can be assigned to (the direct sum of a number of copies of) any vector bundle over $X$ inductively from a $G$-CW skeleton to the one of higher dimension by gluing vector bundles over equivariant cells $G/K\times D^n$, which admit $G$-equivariant structures by the equivariant formality of $G/K$. The clutching map for gluing vector bundles over equivariant cells with the vector bundle over the skeleton can be shown to be homotopy equivalent to a $G$-equivariant one by the same token.
	\end{enumerate}
\end{example}
\subsection{Equivariant formality of homogeneous spaces} Now let us turn to the problem of determining if the action on a general homogeneous space $G/K$, where both $G$ and $K$ are compact and connected, by the left translation by $K$ is equivariantly formal. We call such an action an \emph{isotropy action}. When $K$ is the identity subgroup, equivariant formality of the isotropy action is immediate. On the other hand, the case when $K$ is of the maximal rank is taken care of by Example \ref{fokex} (\ref{eqrank}) and the fact that restricted equivariantly formal action is still equivariantly formal. More subtle is when $K$ sits between the two extreme cases. This problem was taken up in \cite{Sh, ST, Go, GN, Ca}, where the authors either presented a partial characterization of equivariant formality of isotropy actions in cohomological terms or showed that some special examples, e.g. generalized symmetric spaces and when $K$ is a circle, possess such a property by a dimension count argument using Condition (\ref{dimcount}) at the beginning of this section. 

We observe from the previous results that if the isotropy action is equivariantly formal, then $G/K$ itself is formal in the sense of rational homotopy theory (\cite[Theorem A]{CF}). For a formal homogeneous space, its (rational) ordinary $K$-theory ring is generated by the `even' $K$-theory classes represented by homogeneous vector bundles and some special `odd' $K$-theory classes coming from degree $-1$ piece (\cite[Theorem 6.2]{CF}). We have seen from Example \ref{fokex} (\ref{eqrank}) that any homogeneous vector bundle has a natural equivariant structure. However, among the odd $K$-theory generators, only those represented by vector bundles over the suspension of $G/K$ constructed in the same manner as those in Example \ref{fokex} (\ref{oddgen}), but with clutching function $\rho_1\rho_2^{-1}\colon G/K\to U(n)$, $gK\mapsto \rho_1(g)\rho_2(g)^{-1}$, where $\rho_1$ and $\rho_2$ are $G$-representations such that $\rho_1-\rho_2\in\text{ker}(i^*)$, have equivariant structures. By Theorem \ref{RKEF}, determining if the isotropy action is equivariantly formal boils down to deciding if all odd $K$-theory generators are represented by those vector bundles. The latter condition then translates to the following representation theoretic characterization of equivariant formality of isotropy actions which has a different flavor from the previous results. %Recall that the augmentation ideal $I(G)$ of $R(G)$ is the kernel of the augmentation map $R(G)\to \mathbb{Z}$ which takes a representation to its dimension.
\begin{theorem}\textnormal{(cf. \cite[Theorem B]{CF})}. \label{thb}Let $G$ be a connected compact Lie group and $K$ a connected Lie subgroup which acts on $G/K$ by left translation. Let $i^*\colon R(G)\otimes\mathbb{Q}\to R(K)\otimes\mathbb{Q}$ be the restriction map of the rationalized representation ring of $G$ to that of $K$. Then $G/K$ is an equivariantly formal $K$-space if and only if $\text{im }i^*$ is regular at the augmentation ideal $(I(K)\otimes\mathbb{Q})\cap\text{im }i^*$.\footnote{Recall that a Noetherian local ring $R$ with the maximal ideal $\mathfrak{m}$ is said to be regular if its Krull dimension is equal to $\text{dim}_{R/\mathfrak{m}}\mathfrak{m}/\mathfrak{m}^2$. A Noetherian ring $R$ is said to be regular at a prime ideal $\mathfrak{p}$ if the local ring $R_\mathfrak{p}$ is regular.}
\end{theorem}
Theorem \ref{thb} furnishes an elegant and uniform way of showing equivariant formality of some examples of isotropy actions. In the next example, we shall illustrate the utility of this theorem in proving equivariant formality of isotropy actions on generalized homogeneous spaces, which was studied in \cite{Go, GN} in detail on a case-by-case basis using the classification of such spaces. 
\begin{example}(\cite[Proposition 7.9, Proof of Example 2.3 (iii)]{CF}).\label{gensym}
	Let $G$ be a compact, connected, simple and simply-connected Lie group, $\sigma$ an autormorphism of $G$ induced by a Dynkin diagram automorphism, and $K$ the identity component of the subgroup fixed by the automorphism group generated by $\sigma$. The homogeneous space $G/K$ is called a generalized symmetric space. Examples include $SU(n)/SO(n)$ (where $\sigma$ is complex conjugation) and $\text{Spin}(4)/G_2$ (where $\sigma$ is the triality automorphism). Let $\rho_1, \cdots, \rho_n$ be the fundamental representations of $G$ and $\overline{\rho}_i:=\rho_i-\text{dim }\rho_i$ the reduced ones, which are in $I(G)$. The (rationalized) representation ring of $G$ is a polynomial ring generated by $\overline{\rho}_i$, $1\leq i\leq n$ (\cite[Prop. 11.1]{Ho2}). In this sitution the regularity condition in Theorem \ref{thb} amounts to $\text{ker}(i^*)$ being generated by polynomials in $\overline{\rho}_i$ all having nonzero linear terms, which is indeed the case. The group $\langle\sigma\rangle$ acts on the set of reduced fundamental representations and $\text{ker}(i^*)$ is generated by the linear polynomials $\overline{\rho}_i-\overline{\rho}_j$ where $\overline{\rho}_i$ and $\overline{\rho}_j$ are in the same orbit.
\end{example}	
	It is worth mentioning at this point that the notion of affinization discussed immediately after Theorem \ref{specsupp} allows for a nice algebro-geometric interpretation of Theorem \ref{RKEF}. More precisely, the regularity condition means that the map $\text{Spec }R(K)\otimes\mathbb{Q}\to\text{Spec }R(G)\otimes \mathbb{Q}$, which is the affinization of the embedding $K\hookrightarrow G$, is smooth at $I(G)\otimes\mathbb{Q}$. Roughly speaking, on the level of Lie groups, the embedding of $K$ into $G$ is in some sense `smooth' at the identity element, which is the support of $I(G)\otimes\mathbb{Q}$. This `smoothness' is tied with how symmetric the embedding $K\hookrightarrow G$ at the identity element (or equivalently the embedding of the Lie algebra $\mathfrak{k}\hookrightarrow \mathfrak{g}$) is (see Example \ref{gensym} and the remark after \cite[Theorem B]{CF}).

\section{Further developments}
In this last section we would like to give a brief account of recent developments in equivariant twisted $K$-theory. 

Twisted $K$-theory, first introduced in \cite{DP}, is $K$-theory with local coefficient systems called twists, which can be seen as data twisting the transition functions of vector bundles, realized by such models as projective Hilbert space bundles (\cite{AtSe2}) and bundle gerbes (\cite{BCMMS}), and classified by the third integral cohomology group. The introduction of twists are necessary to make up for the absence of $K$-orientation when formulating the Thom isomorphism and Poincar\'e duality. Twisted $K$-theory has found applications in mathematical physics. In string theory, twists play the role of background $B$-fields, and twisted $K$-theory classifies $D$-brane charges and Ramond-Ramond field strength (\cite{W2}). The intriguing phenomenon of $T$-duality also manifests itself in the context of twisted $K$-theory (\cite{BEM}).

Equivariant twisted $K$-theory first received a rigorous treatment in \cite{AtSe2}. There it is defined as the group of homotopy classes of equivariant sections of the bundle associated to the twist whose fibers are the space of Fredholm operators on a Hilbert space, in analogy to the fact that the space of Fredholm operators is a representing space for ordinary $K$-theory. It is also shown that equivariant twists, realized by equivariant projective Hilbert space bundles, are classified by the third equivariant cohomology group. This is a nontrivial result as there is a subtle difference between equivariant bundles and ordinary bundles over the Borel homotopy quotient.

In a series of papers (\cite{FHT, FHT1, FHT2, FHT3}), Freed, Hopkins and Teleman put forth their celebrated theorem (FHT) which identifies the equivariant twisted $K$-theory of a compact Lie group $G$ with conjugation action by itself with its Verlinde algebra, the representation group of its loop group. In these papers they showcased a number of approaches to proving their theorem. Among them is the localization theorem for equivariant twisted $K$-theory. Using this they showed that the sheaf of modules associated with the equivariant twisted $K$-theory of $G$ over $\text{Spec }R(G)\otimes\mathbb{C}$ is a collection of skyscraper sheaves whose support corresponds to the irreducible positive energy representations of the loop group $LG$. They also introduced a generalization of equivariant $K$-theory, namely $K$-theory of groupoids, which include $G$-spaces (called global quotient groupoids) as an example. One can speak of local equivalence of groupoids, and an advantage of working with $K$-theory of groupoids is that it respects local equivalence. The Freed-Hopkins-Teleman Theorem can be heuristically understood through the local equivalence of the global quotient groupoid of $G$ by the conjugation action with the groupoid of the universal proper action by $LG$. The crux of FHT is that it gives an algebro-topological interpretation of the fusion product of the Verlinde algebra, which is defined classically in terms of conformal field theory and algebraic geometry of moduli spaces of Riemann surfaces, through the Pontryagin product of the twisted $K$-theory. Extending this philosophy further, one can, with FHT, adopt the novel point of view that the Verlinde formula, which computes the dimensions of spaces of generalized theta functions for general Riemann surfaces, is an equivariant topological index in twisted $K$-theory (\cite[\S8]{FHT}). This in turn provides a framework for the formulation of the geometric quantization of any quasi-Hamiltonian manifold as the pushforward on equivariant twisted $K$-theory induced by the group-valued moment map, without going through the formalism of Dirac operators on infinite dimensional $LG$-spaces (\cite{Me2}).
\bibliographystyle{amsalpha}

\end{document}